\documentclass{article}[12pt]
\usepackage{amsmath,amssymb,amsthm}
\usepackage{graphicx}

\newenvironment{Proof}{{\it Proof.}}{{$\square$} \vskip.5cm}
\newtheorem{Thm}{Theorem}[section]
\newtheorem{Lem}[Thm]{Lemma}
\newtheorem{Pro}[Thm]{Proposition}
\newtheorem{Cor}[Thm]{Corollary}
\theoremstyle{definition}

\newtheorem{Def}[Thm]{Definition}

\theoremstyle{remark}

\newcommand{\R}{\mathbb{R}}

\newcommand{\cC}{\mathcal{C}}

\newcommand{\cL}{\mathcal{L}}

\newcommand{\cU}{\mathcal{U}}

\newcommand{\ga}{\gamma}

\newcommand{\dd}{\operatorname{\,d}}

\makeatletter
\renewcommand{\subsection}{\@startsection{subsection}{2}%
{0pt}{1ex plus 1ex minus .2ex}
{-0ex}{\normalsize\bf}}
\makeatother
\begin{document}
\title{New lower bound for the distortion\\ of a knotted curve}
\author{T. Bereznyak, P. Svetlov\thanks{The second author was partially supported by
RFBR grant 05-01-00939a.}}
%\date{{\normalsize St.Petersburg branch
%of Steklov Mathematical Institute}}

\maketitle
\begin{abstract} We prove that  distortion of a knotted
curve in $\R^3$ is great than 4.76. This improves a result
obtained by John M. Sullivan and Elizabeth Denne in \cite{DS}.
\end{abstract}
%\twocolumn
\section{The definition and results}
M. Gromov introduced the notion of distortion for submanifolds of
$\R^N$ in \cite{G}, p.113-114. In the case of a simple closed
curve $K$ (closed 1-dimensional submanifold of $\R^3$) distortion
is
$$
\cU(K)=\max\limits_{P,Q\in K}\frac{|PQ|_K}{|PQ|},
$$
where $|PQ|$ is length of the secant $PQ$, and $|PQ|_K$ is length
of the shortest arc in $K$ joining $P$ and $Q$ .

M. Gromov showed that any closed curve has distortion at least
$\pi/2$.  John M. Sullivan and Elizabeth Denne showed in \cite{DS}
that each  knotted curve (=a curve represented a nontrivial knot)
has distortion $>3.99$. (The reader can use the paper \cite{DS} as
an introduction to the topic.) To obtain this bound they use a
special set of ``essential secants'' (see Definition
\ref{ess_sec}). Using similar arguments we obtain the following

\vspace{0.3cm}

{\bf Theorem \ref{main}} {\it Each knotted curve has distortion
$>4.76$.}

\section{Background}
Let $K\subset\R^3$ be an oriented closed 1-dimensional
submanifold. By $\ga_{PQ}$ we denote the arc in $K$ joining points
$P,Q\in K$ following the orientation of $K$.

\begin{Def}{(\cite{DDS})}
Suppose a secant $PQ$ has no interior intersections with the curve
$K$. We say that $\ga_{PQ}$ is an {\it essential arc} if the
closed curve $\ga_{PQ}\cup PQ$ bounds no disk whose interior is
disjoint from $K$.

Suppose a secant $PQ$ has interior intersection with the curve
$K$. We say that $\ga_{PQ}$ is {\it essential} if for any
$\varepsilon>0$ there exist points $P',Q'\subset K$ such that
$|PP'|,|QQ'|<\varepsilon$, the secant $P'Q'$ has no interior
intersections with the curve $K$, and the arc $\ga_{P'Q'}$ is
essential.
\end{Def}

\begin{Def}{(\cite{DDS})}\label{ess_sec}
Given $P,Q\in K$, we say the secant $PQ$ is {\it essential} if
both arcs $\ga_{PQ}$ and $\ga_{QP}$ are essential.
\end{Def}

\begin{Def}{(\cite{DDS})}
We say that an essential arc $\ga_{PQ}$ is {\it
boundary-essential} if there are inessential arcs with endpoints
arbitrarily closed to $P$ and $Q$.
\end{Def}

\begin{Lem}{(\cite{DDS}, Theorem 7.1)}\label{DDS}
If an arc $\ga_{PQ}\subset K$ is boundary-essential then $K$ must
intersect the interior of $PQ$ at some point $O\in\ga_{QP}$ for
which the secants $OP$ and $OQ$ are both essential. $\square$
\end{Lem}

\begin{Lem}\label{cont_ess}
Each essential arc $\ga_{PQ}\subset K$ contains boundary-essential
arcs $\ga_{PS}$ and $\ga_{S'Q}$.
\end{Lem}

\begin{Proof}
Let $s:[0,1]\to \R^3$ be a continuous parametrization of the arc
$\ga_{PQ}$, $s(0)=P$, $s(1)=Q$. There exists $t_0\in (0,1)$ such
that $\ga_{s(t)Q}$ is inessential for any $t>t_0$.
\end{Proof}

\section{$\cU\ge\pi+1$}
Consider a knotted curve $K\subset\R^3$ with distortion $\cU$.
Suppose $L$ is length of a shortest boundary-essential arc in $K$
and $l$ is length of a shortest secant of boundary-essential arcs
in $K$:
$$
L=\min\{|\ga_{AB}|: \ga_{AB}\subset K\quad \mbox{is
boundary-essential}\},
$$
$$
l=\min\{|{AB}|: \ga_{AB}\subset K\quad \mbox{is
boundary-essential}\},
$$
where by $|\ga_{AB}|$ we denote length of the arc $\ga_{AB}$.

%Let  $\ga_{PQ}$ be a boundary-essential arc. By theorem \cite{DDS}
%there is a point $O\in K$ in the interior of $PQ$.
\begin{Thm}\label{avoiding_ball}
Each boundary-essential arc $\ga_{PQ}$ avoids an open ball
centered at some point of the secant $PQ$ of radius
$$
r=\frac{L}{\cU-1}\,.
$$
\end{Thm}
\begin{Proof}
Consider a boundary-essential arc $\ga_{PQ}$. By Lemma \ref{DDS}
the secant $PQ$ intersects $K$ at some point $O$ and secants $PO$,
$OQ$ are both essential. Consider a point $X\in\ga_{PQ}$ and
suppose that $|\ga_{XO}|\le |\ga_{OX}|$. By lemma \ref{cont_ess},
there is a point $S\in \ga_{QO}$ such that the arc $\ga_{SO}$ is
boundary-essential (fig. \ref{fig1}). Note that $|\ga_{SO}|\ge L$,
$|SO|\ge l$ by definitions.
\begin{center}
\begin{figure}
    \includegraphics{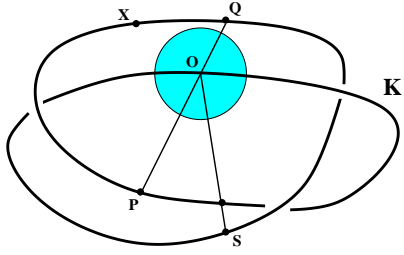}
  \caption{The arc $\ga_{PQ}$ lies outside the ball centered at $O$ of radius
$r=L/(\cU-1)$.}\label{fig1}
\end{figure}
\end{center}
Using obvious inequalities $|\ga_{XS}|\ge |XS|\ge |SO|-|XO|\ge
l-|XO|$, we conclude
$$
\cU\cdot |XO|\ge|\ga_{XO}|=|\ga_{XS}|+|\ga_{SO}|\ge l-|XO|+L.
$$
 So, we have
 \begin{equation}\label{A}
  |XO|\ge \frac{l+L}{1+\cU}\,,
 \end{equation}
 Note, that the inequality (\ref{A}) is true for any point $X\in\ga_{PQ}$.
Applying it to $P,Q\in\ga_{PQ}$ we obtain
\begin{equation}\label{B}
|PQ|=|PO|+|QO|\ge 2\,\frac{l+L}{1+\cU}\,.
\end{equation}
The inequality (\ref{B}) is true for any boundary-essential arc
$\ga_{PQ}$, particulary for boundary-essential arcs with shortest
secants:
$$
l\ge 2\,\frac{l+L}{1+\cU}\,.
$$
So, we have
$$
l\ge\frac{2L}{\cU-1}\,.
$$
The inequality (\ref{A}) is now rewrites as $|XO|\ge L/(\cU-1)$.
\end{Proof}

\begin{Cor}
Distortion of a knotted curve is not less than $\pi+1$.
\end{Cor}
\begin{Proof}
Consider a boundary-essential arc with length $L$. Due to the
previous theorem, we have
$$
L\ge\frac{\pi L}{\cU-1}.
$$
This yields  $\cU\ge\pi+1$.
\end{Proof}

\section{$\cU> 4.76$}
Let $K$ be a knotted curve with distortion $\cU$. In what follows
we assume that $\ga_{PQ}\subset K$ is a shortest
boundary-essential arc. For convenience, rescale the curve so that
$|\ga_{PQ}|=L=\cU-1$. It means that  the radius of the avoiding
ball from theorem \ref{avoiding_ball} is 1. Also we fix a point
$O\in \ga_{QP}$ which lies  in the interior of the secant $PQ$.
\begin{Lem}\label{eval}
\begin{enumerate}
\item $ f\left(|OP|\right)+f\left(|OQ|\right)\le \cU-1, $ where
$f(t)=\sqrt{t^2-1}+\arcsin(1/t)$. \item $1\le |OP|,|OQ|\le
f^{-1}(\cU-1-\pi/2)<\cU-1-\pi/2$.
\end{enumerate}
\end{Lem}

\begin{Proof} The first proposition is an elementary corollary from theorem
\ref{avoiding_ball} (and follows from \cite{DDS}, Lem. 4.3).
 The second one is also elementary:
 note, that $f(t)>t$ if $t\ge 1$ and $|OQ|\ge 1$ by theorem
\ref{avoiding_ball}; we have
$$
f(|OP|)+\frac{\pi}{2}\le f\left(|OP|\right)+f\left(|OQ|\right)\le
\cU-1.
$$
The same arguments are valid for $|OQ|$.
\end{Proof}

Let $P\in\R^3\setminus\{O\}$, $A>0$, $\psi\in [0,\pi/2)$. Consider
a family $\cC(P,A,\psi)$ of arc-length parameterized curves
$C:[0,a_C]\to\R^3$ such that $C(0)$ lies on the axis $OP$ and
\begin{equation}\label{bsc}
|OC(t)|\ge A+t\cdot\sin{\psi}.
\end{equation}
\begin{Lem}\label{body}
Each curve $C\in\cC(P,A,\psi)$ avoids the body
$$
B=
\left\{X\in\R^3:|OX|<A\cdot e^{\angle{POX}\tan\psi}
\right\}.
$$
\end{Lem}
\begin{Proof}
We use spherical coordinate system gave by a splitting
$$\R^3\setminus \{O\}=(0,+\infty)\times S^2.
$$
Each point $X\in\R^3\setminus \{O\}$ has two coordinates, say
$(\rho_X,\omega_X)$. The first one  $\rho_X$ is the distance
$|OX|$, and the second one $\omega_X$ is the projection of $X$ on
the unit sphere centered at $O$. It is obvious that
$$
\cC(P,A,\psi)=\Bigl\{
\{\rho(t),\omega(t)\}:\,\omega(0)=\omega_P,\, \rho(t)\ge
A+t\sin\psi\Bigr\}.
$$
Each curve from the boundary of $\cC(P,A,\psi)$ satisfies
$$
\left\{
\begin{array}{ll}
\rho(0)=A,\,\omega(0)=\omega_P&\mbox{(the initial conditions)}\\
\rho'(t)=\sin\psi&\mbox{(differentiation of (\ref{bsc}))}\\
\rho'^2+\rho^2\omega'^2=1&\mbox{(arc-lenght parametrization)}
\end{array}\right.
$$
Two last equations give $(\ln\rho)'=\omega'\tan\psi$ that yields
$$
\ln\rho_X-\ln A=\tan\psi\int^{X}_{P}\dd \omega\ge \tan\psi\angle POX
$$

\end{Proof}

By definition, put
\begin{equation}\label{C}
x_Q=\frac{\cU+1-|OQ|}{\cU},\quad x_P=\frac{\cU+1-|OP|}{\cU},
\end{equation}
and note that $x_P, x_Q\in(0,1]$ (lemma \ref{eval}). Consider
rotation bodies
$$
B_P=\left\{X\in\R^3:\,|OX|<x_P\cdot exp\left(\frac{\angle
POX}{\sqrt{\cU^2-1}} \right)\right\},
$$
$$
B_Q=\left\{X\in\R^3:\,|OX|<x_Q\cdot exp\left(\frac{\angle
QOX}{\sqrt{\cU^2-1}} \right)\right\}.
$$

\begin{center}
\begin{figure}
    \includegraphics{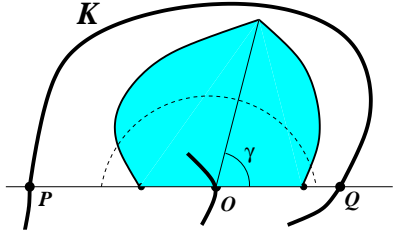}
  \caption{ The avoiding body. Here  $\gamma={\pi}/{2}-\sqrt{\cU^2-1}\cdot\ln \sqrt{x_Q/x_P}$. The punctured line shows the avoiding ball
  from lemma \ref{avoiding_ball}.}\label{fig2}
\end{figure}
\end{center}

\begin{Pro}\label{avoiding_body}
Each boundary-essential arc $\ga_{PQ}$ avoids the intersection $B_{P}\cap B_Q$ (see fig.\ref{fig2}).
\end{Pro}
\begin{Proof}
Let $X\in\ga_{PQ}$ and $|\ga_{XO}|\le |\ga_{OX}|$. We have
$\cU\cdot |XO|\ge |\ga_{XO}|=$
$$
=|\ga_{XQ}|+|\ga_{QS}|+|\ga_{SO}|\ge
|\ga_{XQ}|+|OS|-|OQ|+L\ge $$$$\ge L+\frac{2L}{\cU-1}-|OQ|+|\ga_{XQ}|=
\cU+1-|OQ|+|\ga_{XQ}|.
$$
It means that
\begin{equation}
|OX|\ge  \left\{\begin{array}{lll}
x_Q+|\ga_{XQ}|/\cU,&\mbox{if}& |\ga_{XO}|\le |\ga_{OX}|\\
&&\\
x_P+|\ga_{PX}|/\cU,&\mbox{if}& |\ga_{XO}|\ge |\ga_{OX}|                \end{array}\right..
\end{equation}
Applying the lemma \ref{body} we are done.
\end{Proof}

From now on we assume that $\cU\le 5$.

 Now we compute length of
the shortest curve joining points $P$ and $Q$ outside the avoiding
body from proposition \ref{avoiding_body}. Obviously, the curve is
a plane one.

\begin{center}
\begin{figure}
    \includegraphics{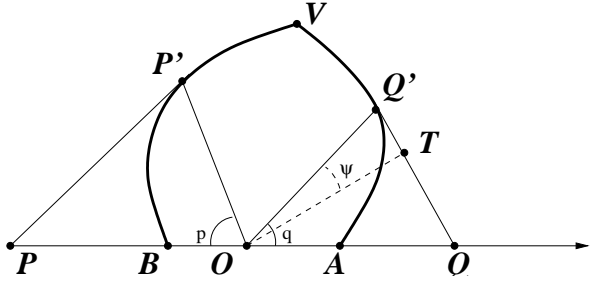}
  \caption{The shortest curve $PP'VQ'Q$ joining points $P$ and $Q$
  outside the avoiding body $BP'VQ'A$. Here $|OA|=x_Q$, and $|OB|=x_P$}\label{fig3}
\end{figure}
\end{center}

\begin{Lem} 1. The shortest curve joining points $P$ and $Q$ outside the avoiding
body from proposition \ref{avoiding_body} consists on four
nonempty arcs $\ga_{PP'}$, $\ga_{P'V}$, $\ga_{VQ'}$, and
$\ga_{Q'Q}$ (see fig. \ref{fig3}).

2. Length of the curve is
$$
F(Q,P;\cU)=\frac{2e^{\frac{\pi}{2}\tan\psi}}{\sin{\psi}}\sqrt{x_Px_Q}-\cot\psi\left(
\frac{x_Qe^{q\tan\psi}\cos
q}{\cos(q-\psi)}+\frac{x_Pe^{p\tan\psi}\cos p}{\cos(p-\psi)}
\right)\,,
$$

where $\psi=\arcsin(1/\cU)$, $p=\angle(POP')$, $q=\angle(QOQ')$.
\end{Lem}
\begin{Proof}
It can easily be checked that $\angle(TOQ')=\psi$,
$|OQ'|=x_Qe^{q\tan\psi}$, and
$$
|OQ|=g_\psi(q)\cdot x_Q,\quad\mbox{where}\quad g_\psi(q)=\frac{e^{q\cdot\tan{\psi}}\cdot\cos{\psi}}{\cos(\psi-q)}.
$$
Taking into account (\ref{C}) we obtain the following dependence
between $x_Q$ and $q$:
\begin{equation}\label{D}
x_Q=\frac{\cU+1}{\cU+g_\psi(q)}\,.
\end{equation}
Note, that
$$
\frac{\dd g_\psi(q)}{\dd q}=\frac{e^{q\cdot\tan{\psi}}\sin{q}}{\cos^2(\psi-q)}.
$$
 It means that $x_Q$ is a decreasing function on $q\in[0,\psi+\pi/2)$.

Now we claim that if $\cU\le 5$, then $q<\gamma=\angle(QOV)<\pi-p$
(see Appendix).  It proves the first part of the lemma.

The remain part is some computation:
$$
|QQ'|=\sqrt{|QO|^2+|Q'O|^2-2|QO|\cdot|Q'O|\cos{q}}
%=x_Qe^{q\tan\psi}\sqrt{\left(\frac{\cos\psi}{\cos(q-\psi)}\right)^2+1-2\frac{\cos\psi\cos{q}}{\cos(q-\psi)}}=\\
=x_Q\frac{e^{q\cdot\tan{\psi}}\sin{q}}{\cos(q-\psi)},
$$
$$
|PP'|=\sqrt{|PO|^2+|P'O|^2-2|PO|\cdot|P'O|\cos{p}}
%=x_Qe^{q\tan\psi}\sqrt{\left(\frac{\cos\psi}{\cos(q-\psi)}\right)^2+1-2\frac{\cos\psi\cos{q}}{\cos(q-\psi)}}=\\
=x_P\frac{e^{p\cdot\tan{\psi}}\sin{p}}{\cos(p-\psi)},
$$
$$
|\ga_{Q'V}|=\cU
x_Q(e^{\gamma\cdot\tan\psi}-e^{q\cdot\tan\psi})=\cU\sqrt{x_Px_Q}e^{\frac{\pi}{2}\tan\psi}-
\cU x_Qe^{q\cdot\tan\psi},
$$
$$
|\ga_{VP'}|=\cU
x_P(e^{(\pi-\gamma)\cdot\tan\psi}-e^{p\cdot\tan\psi})=\cU\sqrt{x_Px_Q}e^{\frac{\pi}{2}\tan\psi}-
\cU x_Pe^{p\cdot\tan\psi}.
$$

\end{Proof}
\begin{Thm}\label{main}
Each knotted curve has distortion $>4.76$.
\end{Thm}

\begin{Proof}
Consider the function $\cL:[\pi+1,5]\to\R$
$$
\cL(t)=\min\left\{ F(P,Q;t): f(|OP|)+f(|OQ|)\le t-1 \right\}.
$$
It  follows easily that
$$
\cU\ge \max\{t\in[\pi+1,5]: \cL(t)\ge t-1\}.
$$
Using mathematical software we found that $\max\{t\in[\pi+1,5]:
\cL(t)\ge t-1\}> 4.76$.
\end{Proof}

\section{Appendix}
\begin{Lem} If $\cU\le 5$ then
$x_P,x_Q\in (0.76,1]$,
and
 $p,q\in[0,1.32)$.
\end{Lem}
\begin{Proof} As it follows from lemma \ref{eval},
variables $x_P,x_Q$ lie in the domain
$$
\left\{x_Q,x_P\le 1, f(\cU+1-\cU x_Q)+ f(\cU+1-\cU
x_P)\le\cU-1\right\}
$$
 Let
$u(t)=f^{-1}(t-1-\pi/2)$, $t\ge 1$. Consider the function
$$
x(t)=1-\frac{u(t)-1}{t}.
$$
It is easily shown that $x'(t)<0$. Therefore
$$
0.76<x(5)\le x(\cU)\le x_Q,x_P\le 1.
$$
By (\ref{D}) we have
$$
\frac{\cU+1}{\cU+g_\psi(q)}\ge x(\cU).
$$
This inequality implies $0\le q< 1.32$ whenever $\cU\le 5$.
\end{Proof}

The desired inequality $\ga>q$ is equivalent to the following
$$
\pi/2+\frac{\cot\psi}{2}\ln x_P>q+\frac{\cot\psi}{2}\ln x_Q.
$$
The left hand side is not less than
$$
\pi/2+\frac{\cot\psi}{2}\ln x(\cU)\ge \pi/2+\sqrt{6}\ln x(5)>0.9
$$
\begin{Lem} If $\cU\le 5$ then
$q+\frac{\cot\psi}{2}\ln x_Q<0.9$.
\end{Lem}
\begin{Proof}
It can easily be checked that the negative function
$d(q,U)=\frac{\cot\psi}{2}\ln x_Q$ is a decreasing function on
both $q\in [0.9,1.32]$ and $\cU\in [\pi+1,5]$. The following
sequence completes the proof:
$$
\begin{array}{lll} q+\frac{\cot\psi}{2}\ln x_Q<0.9&\mbox{if}&
q\in[0,0.9)\\
&&\\ q+\frac{\cot\psi}{2}\ln x_Q\le
q+d(0.9,\pi+1)<q-0.19<0.9&\mbox{if}& q\in[0.9,1.09)\\
&&\\
q+\frac{\cot\psi}{2}\ln x_Q\le
q+d(1.09,\pi+1)<q-0.33<0.9&\mbox{if}& q\in[0.9,1.23)\\
&&\\
q+\frac{\cot\psi}{2}\ln x_Q\le
q+d(1.23,\pi+1)<q-0.48<0.9&\mbox{if}& q\in[1.23,1.32]
\end{array}.
$$

\end{Proof}

%%%%%%%%%%%%%%%%%%%%%%%%%%%%%%%%%%%%%%%%%%%%%%%%%%%%%%%%%%%%%%%%%%

\bigskip

\noindent Math. Dep. of St.-Petersburg State University

\noindent Staryi Petergof, Universitetsky pr., 28

\noindent 198504, St.-Petersburg, Russia

\noindent {\tt tarasber@yandex.ru}

\bigskip

\noindent St.-Petersburg Dept. of Steklov Math. Institute

\noindent Fontanka 27,

\noindent 191011, St.-Petersburg, Russia

\noindent {\tt svetlov@pdmi.ras.ru}

\end{document}